\documentclass{amsart}

\newtheorem{teo}{Theorem}[section]
\newtheorem{prop}[teo]{Proposition}
\newtheorem{cor}[teo]{Corollary}
\newtheorem{lema}[teo]{Lemma}
\newtheorem{claim}[teo]{Claim}

\theoremstyle{definition}
\newtheorem{dfn}[teo]{Definition}
\newtheorem{ex}[teo]{Example}

\theoremstyle{remark}
\newtheorem{rem}[teo]{Remark}

\numberwithin{equation}{section}

\newcommand{\length}{\ensuremath{\mathrm{length}}}
\newcommand{\WS}{\ensuremath{W_{\Sigma}}}
\newcommand{\WSt}{\ensuremath{W_{\widetilde{\Sigma}}}}
\newcommand{\WSr}{\ensuremath{W_{\Sigma_{r}}}}
\newcommand{\WSrt}{\ensuremath{W_{\widetilde{\Sigma}_{r}}}}

\newcommand{\s}{\ensuremath{\varsigma}}
\newcommand{\interior}{\ensuremath{\mathrm{int} }}
\newcommand{\dist}{\ensuremath{\mathrm{dist} }}

\newcommand{\tub}{\ensuremath{\mathrm{Tub} }}

\newcommand{\F}{\ensuremath{\mathcal{F}}}
\newcommand{\singularF}{\ensuremath{\mathcal{X}_{F}}}

\newcommand{\dank}{\textsf{Acknowledgments:\ }} 

\newcommand{\eps}{\varepsilon}

\begin{document}

\title{Singular riemannian foliations on simply connected spaces}

\author{Marcos M. Alexandrino \and Dirk T\"oben} 
 
\thanks{ M.M. Alexandrino was supported in part by DFG-Schwerpunkt Globale Differentialgeometrie\\ D.T\"{o}ben was supported in part by PRONEX-RIO}

\address{Marcos M. Alexandrino\\Departamento de Matem\'{a}tica, PUC-Rio, 
 Rua Marqu\^{e}s de S\~{a}o Vicente, 225, 
  22453-900, Rio de Janeiro, Brazil}
\email{malex@mat.puc-rio.br}

\address{Dirk T\"oben\\
Mathematisches Institut, Universit\"at zu K\"oln, Weyertal 86-90, 50931 K\"oln, Germany}
\email{dtoeben@math.uni-koeln.de}


\date{November, 2004.}


\begin{abstract}
A singular foliation on a complete riemannian manifold is said to be riemannian if each geodesic that is perpendicular at one point to a leaf remains perpendicular to every leaf it meets. The singular foliation is said to admit sections if 
each regular point is contained in a totally geodesic complete immersed submanifold that meets every leaf orthogonally and whose dimension is the codimension of the regular leaves.
A typical example of such a singular foliation is the partition by orbits of a polar action, e.g. the orbits of the adjoint action of a compact Lie group on itself.

We prove that a singular riemannian foliation with compact leaves  that admits sections on a simply connected space has no exceptional leaves, i.e., each  regular leaf  has trivial normal holonomy. We also prove that there exists a convex fundamental domain in each section of the foliation and in particular that the space of leaves is a convex Coxeter orbifold.
\end{abstract}

\maketitle

\section{Introduction}

In this section we will recall the concept of a singular riemannian foliation with sections, review typical examples
and  state our main results as Theorem \ref{teo-excepcional-leaf}
 and Theorem \ref{teo-fundamental-domain}.

We start by recalling the definition of a singular riemannian foliation (see  the  book of  P. Molino \cite{Molino}).

\begin{dfn}
 A partition $\F$ of a complete riemannian manifold $M$ by connected immersed submanifolds (the \emph{leaves}) is called a \emph{singular foliation} of $M$ if it verifies condition (1) and {\it singular riemannian foliation} if it verifies condition (1) and (2):

\begin{enumerate}
\item $\F$ is \emph{singular},
i.e., the module $\singularF$ of smooth vector fields on $M$ that are tangent at each point to the corresponding leaf acts transitively on each leaf. In other words, for each leaf $L$ and each $v\in TL$ with footpoint $p,$ there is $X\in \singularF$ with $X(p)=v$.
\item  The partition is \emph{transnormal}, i.e., every geodesic that is perpendicular at one point to a leaf remains perpendicular to every leaf it meets.
\end{enumerate}
\end{dfn}

Let $\F$ be a singular riemannian foliation on a complete riemannian manifold $M.$  A leaf $L$ of $\F$ (and each point in $L$) is called \emph{regular} if the dimension of $L$ is maximal, otherwise $L$ is called {\it singular}. Let $L$ be an immersed submanifold of a riemannian manifold $M.$  A section $\xi$ of  the normal bundle $\nu L$ is said  to be a \emph{parallel normal field} along $L$ if  $\nabla^{\nu}\xi\equiv 0,$ where $\nabla^{\nu}$ is the normal connection.  $L$  is said to have a {\it globally flat normal bundle},  if the holonomy  of the normal bundle  $\nu L$ is trivial, i.e., if any normal vector can be extended to a globally defined parallel normal field. 
Examples of submanifolds with flat normal bundle are the regular leaves of the singular foliation defined below. 

\begin{dfn}[s.r.f.s.]
 Let $\F$ be a singular riemannian foliation on a complete riemannian manifold $M.$
$\F$ is said to be a \emph{singular riemannian foliation with sections} (s.r.f.s. for short) if for each regular point $p,$ the set $\Sigma :=\exp_{p}(\nu_p L_{p})$ is a complete immersed submanifold that meets each leaf orthogonally. $\Sigma$ is called a \emph{section}.
\end{dfn}

Singular riemannian foliations with sections were first studied by the first author in \cite{Alex1}, \cite{Alex2} and \cite{Alex3} and continued to be studied by the second author in  \cite{Toeben1} and \cite{Toeben2}. We will recall some properties of s.r.f.s. in the next section. 

Typical examples of singular riemannian foliations with sections are the set of orbits of a polar action, parallel submanifolds of an isoparametric submanifold in a space form and parallel submanifolds of an equifocal submanifold with flat sections in a simply connected compact symmetric space. We will now briefly recall these notions.
An isometric action of a compact Lie group $G$ on a complete riemannian manifold $M$ is called \emph{polar} if there exists a complete immersed submanifold $\Sigma$ of $M$ 
that meets all $G$-orbits orthogonally and whose dimension is equal to the codimension of the regular $G$-orbits. 
We call $\Sigma$ a \emph{section}. A typical example of a polar  action is a compact Lie group with a biinvariant metric that acts on itself  by conjugation. In this case the maximal tori are the sections.

A submanifold of a real space form is called \emph{isoparametric} if its normal bundle is flat and if the principal curvatures along any parallel normal vector field are constant. In \cite{TTh1} C.L. Terng and G. Thorbergsson introduced the class of equifocal submanifolds in simply connected, compact symmetric spaces that have similar properties as isoparametric ones. In \cite{Toeben2} the second author gave a necessary and sufficient condition for an equifocal submanifold (which he calls submanifold with parallel focal structure) in an arbitrary ambient space to induce a s.r.f.s. by parallel and focal submanifolds. From this he derived similar properties. The history of isoparametric submanifolds and their generalizations can be found in the survey \cite{Th} of G. Thorbergsson (see also \cite{Th2}).

\begin{dfn}
\label{dfn-equifocal}
A connected immersed submanifold  $L$  of a complete riemannian manifold $M$ is called \emph{equifocal} if the conditions below are satisfied.
\begin{enumerate}
\item The normal bundle $\nu L$ is  flat.
\item Let $U\subset L$ be a neighborhood small enough such that $\nu(U)$ is globally flat and let $\xi$ be a parallel normal field on $U$. Then the derivative of the map $\eta_{\xi}:U\rightarrow M,$ defined by $\eta_{\xi}(x):=\exp_{x}(\xi),$ has constant rank.
\item $L$ has sections, i.e., for all $p\in L$ there exists a complete totally geodesic immersed submanifold $\Sigma$, the {\it section},  such that $\nu_{p}L=T_{p}\Sigma.$
\end{enumerate}
\end{dfn}
%%%%%%%%pre-introducao

Finally we are ready to state our main results.

\begin{teo}
\label{teo-excepcional-leaf}
Let $\F$ be a s.r.f.s. on a simply connected complete riemannian manifold $M.$ Suppose that the leaves of $\F$ are compact. Then each regular leaf has trivial normal holonomy.
\end{teo}

This result generalizes a theorem of the second author \cite{Toeben1}, who proved the result under the additional condition that the sections are symmetric or do not have any conjugate points. Compare also with Lemma 1A.3 in \cite{PoTh}.

\begin{teo}
\label{teo-fundamental-domain}
Let $\F$ be a s.r.f.s. on a simply connected riemannian manifold $M.$ Suppose also that the leaves of $\F$ are compact. Then
\begin{enumerate}
\item $M/F$  is a simply connected Coxeter orbifold.
\item Let $\Sigma$ be a section of $\F$ and $\Pi:M\rightarrow M/\F$ the canonical projection. Let $\Omega$ be a connected component of the set of regular points in $\Sigma$. Then  $\Pi:\Omega\rightarrow M_{r}/\F$ and  $\Pi:{\overline{\Omega}}\rightarrow M/\F$ are homeomorphisms, where $M_r$ denotes the set of regular points in $M$. In addition $\Omega$ is convex, i.e., for any two points $p$ and $q$ in $\Omega$ every minimal geodesic segment between $p$ and $q$ lies in $\Omega.$
\end{enumerate}
\end{teo}

This paper is organized as follows. In Section 2 we review some facts about s.r.f.s. and fix the notation. In Section 3 we give some results about the space of leaves and prove Theorem \ref{teo-fundamental-domain}. In Section 4 we introduce the concept of a
transverse frame bundle associated to a s.r.f.s. In Section 5 we prove Theorem \ref{teo-excepcional-leaf} using the bundle defined in Section 4. Section 6, the Appendix, is based on E. Salem \cite[Appendix D]{Molino}. There we recall the definitions of pseudogroups, orbifolds, $W$-loops and fundamental groups of $W$-loops.

\dank We are very grateful to  Professor Gudlaugur Thorbergsson for his consistent support  and for helpful discussions. In particular the first author would like to thank Professor Thorbergsson for another invitation to the Mathematisches Institut der Universit\"{a}t zu K\"{o}ln.  The second author would like to thank Professor Ricardo Sa Earp for the invitation to the Departamento de Matem\'{a}tica of PUC-Rio.

\section{Facts about s.r.f.s.}

In this section we recall some results about s.r.f.s. that  will be used in this work. Details can be found in \cite{Alex2} or in \cite{Toeben2}. 
Throughout this section we assume that $\F$ is a s.r.f.s. on a complete riemannian manifold $M.$

Let us start with a result that relates s.r.f.s. to equifocal submanifolds (see Definition \ref{dfn-equifocal} to recall the definitions of equifocal submanifolds and the endpoint map $\eta_{\xi}$).

\begin{teo}
\label{frss-eh-equifocal}
The regular leaves of $\F$ are  equifocal. In particular, the union of regular leaves that 
have trivial normal holonomy is an open and dense set in $M,$ provided that all the leaves are compact.
\end{teo}

A consequence of the first statement of the theorem is that given a regular leaf $L$ with trivial holonomy then we can reconstruct $\F$ by taking all parallel submanifolds of $L.$ More precisely we have

\begin{cor}
\label{cor-map-paralelo}
Let $L$ be a regular leaf of  $\F.$ 
\begin{enumerate}
\item[a)] Let  $\beta$ be a smooth curve of $L$ and  $\xi$ a parallel normal field along  $\beta$. Then the curve $\eta_{\xi}$ belongs to a leaf of \F.
\item[b)] Let  $L$ be a regular leaf with trivial holonomy and $\Xi$ denote the set of all parallel normal fields along $L.$ Then $\F=\{\eta_{\xi}(L)\}_{\xi\in \, \Xi}.$ 
In particular, if $\xi$ is a parallel normal field along $L$ then the endpoint map $\eta_{\xi}:L\rightarrow L_{q}$ is surjective, where $q=\eta_{\xi}(x)$ for $x\in L.$
\end{enumerate}
\end{cor}

Corollary \ref{cor-map-paralelo} allows us to define parallel displacement maps (see definition below) which will be very useful to study $\F.$  

Let $\Sigma$ be a section and $p\in\Sigma.$ Let $B$ a convex normal ball neighborhood of $p$ in $M.$ A connected component of $B\cap \Sigma$ that contains $p$ is called {\it local section} (centered at $p$).

\begin{prop}[Parallel Displacement Maps]
\label{prop-holonomia-singular}
Let $L_{p}$ be a regular leaf,  $\beta$ a smooth curve in $L_{p}$ and let $[\beta]$ denote the homotopy class of $\beta.$ Let $U$ be a local section centered at $p=\beta(0).$ Then there exists a local section $V$ centered at $\beta(1)$ and 
an isometry $\varphi_{[\beta]}:U\rightarrow V$  that has the following properties:
\begin{enumerate}
\item[1)]$\varphi_{[\beta]}(x)\in L_{x}$ for each $x\in U,$
\item[2)]$d\varphi_{[\beta]}\xi(0)=\xi(1),$ where $\xi$ is a  parallel normal field along $\beta.$
\end{enumerate}
\end{prop}

An isometry as above is called \emph{parallel displacement map along $\beta$}. 
 We remark that, in the definition of the parallel displacement map, singular points can be contained in the domain $U.$  
If the domain $U$ and the range $V$ are sufficiently small, then the parallel displacement map coincides with the holonomy map along $\beta.$

Now we  recall some results about the local structure of $\F$, in particular about the structure of the set of singular points in a local section.

Let $p$ be a point of a leaf $L$ of a singular riemannian foliation and let $P$ be a relatively compact, simply connected, open neighborhood of $p$ in $L$. Then there is $\eps>0$ such that $\exp$ restricted to $\nu^\eps P=\{X\in \nu P\mid \|X\|<\eps\}$ is a diffeomorphism onto the tube $\tub (P)$ of radius $\eps$. We call this tube a {\it distinguished} neighborhood of $P$. We write $S_q=\exp(\nu^\eps_qP)$ for $q\in P$. 
We call $S_{q}$ a \emph{slice} through $q.$ 

Note that, if $q$ is a singular point, the restriction $\F|S_q$ of $\F$ to $S_{q}$ is also a singular foliation. In fact, since $\singularF$ acts transitively on the leaves,  the plaques (the connected components of the leaves intersected with $\tub(P)$) are transversal to $S_q.$ 

The relation between a slice $S_{q},$ a local section and $\F|S_q$ are given by the next result, proved by the first author in \cite{Alex2} and \cite{Alex3}. 

\begin{teo}[Slice Theorem]
\label{sliceteorema}
Let $q$ be a singular point of $M$ and $S_{q}$ a slice at $q.$ Then
\begin{enumerate}
\item[a)] Let $\epsilon$ be the radius of the slice $S_{q}.$ Denote  $\Lambda(q)$  the set of 
local sections $\sigma$ containing $q,$ such that $\dist(p,q)<\epsilon$ for each $p\in\sigma.$ 
 Then  $S_{q}= \cup_{\sigma\in\Lambda (q)}\, \sigma.$
%$S_{q}= \cup_{\sigma\in\Lambda (q)}\, \sigma,$ where $\Lambda(q)$  is the set of all local %sections $\sigma$ that contain $q,$ 
\item[b)] $S_{x}\subset S_{q}$ for all $x\in S_{q}.$
\item[c)] $\F|S_q$ is a s.r.f.s. on $S_{q}$ with the induced metric of $M.$
\item[d)] $\F|S_q$ is diffeomorphic  to an isoparametric foliation on an open set of $\mathbf{R}^{n},$  where $n$ is the dimension of   $S_{q}.$
\end{enumerate}
\end{teo}

 From d) it is not difficult to derive the following corollary.

\begin{cor}
\label{estratificacao-singular}
Let $\sigma$ be a local section. Then the set of singular points of $\F$ contained in $\sigma$ is a finite union of totally geodesic hypersurfaces. 
These hypersurfaces are sent by a diffeomorphism to focal hyperplanes 
 contained in a section of an isoparametric foliation on an open set of a euclidean space.   
\end{cor}
We will call the set of singular points of $\F$ contained in $\sigma$ the \emph{singular stratum of the local section} $\sigma$. 
Recall that $M_{r}$ denotes the set of regular points in $M.$ A \emph{Weyl Chamber} of a local section $\sigma$ is the closure in $\sigma$ of a connected component of $M_{r}\cap\sigma.$  

\begin{cor}
\label{convexidadedeWeyl}
A Weyl Chamber of a local section is a convex set.
\end{cor}

\begin{cor}
\label{pontos-singulares-na-geodesica}
Let $\gamma$ be a geodesic orthogonal to a regular leaf. Then the  singular points are isolated on $\gamma.$ 
\end{cor}

\begin{rem}
The above fact can be proved without the Slice Theorem. In fact it was used in \cite{Alex2} and  \cite{Alex3} to prove the Slice Theorem.
\end{rem}

\begin{cor}[Trivialization of $\F$]
\label{transnormalprojection}
For each $q\in M$ there exists a tubular neighborhood $\tub (P_{q})$ of a plaque $P_{q},$ a Weyl chamber $C$ that contains $q$ and a continuous surjective map  $T:\tub (P_{q})\rightarrow  C$ 
with the following property. For each $x\in\tub(P_{q}),$ the point $T(x)$ is the unique point of the intersection of the plaque $P_{x}$ with $C,$ where $P_{x}$ is the plaque in $\tub(P_{q})$ that contains $x.$ 
\end{cor}

\begin{rem}
We can also describe the plaques of $\F$ as level sets of a smooth map $T:\tub(P_{q})\rightarrow \textbf{R}^{k},$ where $k$ is the dimension of the sections. For each regular value $c\in T(\tub(P_{q}))$ there exists a neighborhood $V$ of $T^{-1}(c)$ in $M$ such that $T|_{V}:V\rightarrow T(V)$ is a riemannian submersion with integrable horizontal distribution  for some metric in $T(V).$
These maps are called transnormal maps and more details about them can be found in \cite{Alex1}.
\end{rem}

The Slice Theorem establishes a relation between s.r.f.s. and isoparametric foliations. By analogy with the classical theory of isoparametric submanifolds, it is 
 natural to ask if we can define a (generalized) Weyl group action on $\sigma.$ The next definitions and results answer this question.

\begin{dfn}[Weyl Pseudogroup $W$]
\label{definitionWeylPseudogroup}
 The pseudosubgroup  generated by all parallel displacement maps $\varphi_{[\beta]}$ such that $\beta(0)$ and $\beta(1)$ belong to the same local section $\sigma$ is called \emph{generalized Weyl pseudogroup of} $\sigma.$ Let $W_{\sigma}$ denote this pseudogroup. In a similar way we define $W_{\Sigma}$ for a section $\Sigma.$ Given a  slice $S$ we will define $W_{S}$ as the set of  all parallel displacement maps   $\varphi_{[\beta]}$ such that $\beta$ is contained in the slice $S.$
\end{dfn}

\begin{rem}
To recall the definitions of pseudogroups, orbifolds and $W$-deformation loops see the Appendix. 
\end{rem}

\begin{prop}
Let $\sigma$ be a local section. Then the reflections in the hypersurfaces of the singular stratum of the local section $\sigma$ leave $\F|\sigma$  invariant. Moreover these reflections are elements of $W_{\sigma}.$
\end{prop}

One can construct an example of a s.r.f.s. by suspension such that  $W_{\sigma}$ is larger than the pseudogroup generated by the reflections  in the  hypersurfaces of the singular stratum of $\sigma.$ The next result gives a sufficient condition 
to guarantee that both pseudogroups coincide.

\begin{prop}
\label{Weylgrouparereflection}
Suppose that each leaf of $\F$ is compact and has trivial normal holonomy. Let $\sigma$ be a local section.  Then  $W_{\sigma}$ is generated by the reflections  in the hypersurfaces of the singular stratum of the local section. 
\end{prop}

%%%%%%%%%%%%%%%%%%%%%%%%%%%
%pseudogroup (definition see Molino pag 266, H -loog pag 273 and 274). 

%%%%%%%%%%%%%%%%%%%%%%%%%%%%%%%%%%%%%%%%%%%%%%%%%%%%%%%%%%%%%%%%%%%%%%%%

\section{The space of leaves $M/\F$}
In this section we study the space of leaves $M/F$ and prove Theorem \ref{teo-fundamental-domain}. We will use definitions and results of the last section, as well as the definitions of pseudogroups, orbifolds and $W$-deformation loops, which are explained in the Appendix.

Throughout this section we assume that $\F$ is a s.r.f.s. on a complete riemannian manifold $M$
 and $\Pi: M\rightarrow M/\F$ is the natural projection.

\begin{prop}
\label{orbifoldcoxeter}
Let $\Sigma$ be a section of $\F$ and let $W_{\Sigma}$ be the generalized Weyl pseudogroup of $\Sigma.$ Then $\Sigma/W_{\Sigma} = M/\F.$ If the leaves of $\F$ are compact and the holonomy of each regular leaf is trivial, then $M/F$  is a Coxeter orbifold.
\end{prop}
\begin{proof}
A leaf intersects a given section $\Sigma$ in a $W_\Sigma$-orbit. This defines a continuous, bijective map $H:M/\F\to \Sigma/W_\Sigma$. We want to show that its inverse map is continuous. We consider $W_\Sigma(p)$ for a point $p\in\Sigma$. Let $\sigma$ be a local section centered at $p$ and let $C$ be a Weyl chamber. Clearly $\sigma/W_\sigma$ is a neighborhood of $W_\Sigma(p)$ in $\Sigma/W_\Sigma$. By Proposition \ref{Weylgrouparereflection} it is homeomorphic to $C$ and $\Sigma/W_\Sigma$ is a Coxeter orbifold. This gives a continuous map $\iota_p:\sigma/W_\sigma\to C$. The restriction of the inverse of $H$ to $\sigma/W_\sigma$ is $\Pi\circ \iota_p,$ which  is  continuous. Therefore the inverse of the map $H$ is continuous.
\end{proof}

\begin{prop}
\label{PropHomOfFundamentalGroupofPseudogroup}
Let $\Sigma$ be a section of $\F$ and $p$ a point in $\Sigma.$ Let $W_{\Sigma}$ be the  generalized Weyl pseudogroup of $\Sigma.$ Then there exists a surjective homomorphism $f:\pi_{1}(M,p)\rightarrow \pi_{1}(W_{\Sigma},p).$
\end{prop}
\begin{proof}
The proof is very similar to the proof of E. Salem for the case of regular riemannian foliations (see \cite{Molino}, page 275).

We will first define the homomorphism $f$. According to Corollary \ref{transnormalprojection} there exists an open covering $\{U_{j}\}$ of $M$ such that, for each $j,$ the plaques of $U_{j}$  are preimages of a continuous map $T_{j}:U_{j}\rightarrow C_{j},$ where $C_{j}$ is a Weyl chamber in $U_{j}.$ Let $\alpha$ be a loop in $M$ based at $p\in\Sigma.$  We can choose a subdivision $0=t_{0}<\cdots<t_{n}=1$ such that for each $i$ there exists a $j_{i}$ with $\alpha|_{[t_{i-1},t_{i}]}\subset U_{j_{i}}.$ Define $\tilde{\alpha}_{i}:=T_{j_{i}}\circ \alpha|_{[t_{i-1},t_{i}]}.$ Now for each $\tilde{\alpha}_{i}$ we can find a parallel displacement map $\varphi_{i}$ such that $c_{i}=\varphi_{i}\circ \tilde{\alpha}_{i}$ is a curve in $\Sigma.$ Since $c_{i}(t_{i})$ and $c_{i+1}(t_{i})$ are in the same leaf there exists $w_{i}\in W_{\Sigma}$ such that $w_{i}c_{i}(t_{i})=c_{i+1}(t_{i}).$   
Therefore $(w_{i},c_{i})_{0\leq i\leq n}$ is a $W_{\Sigma}$-loop. The equivalence class of this $W_{\Sigma}$-loop depends only on the loop $\alpha$ and not on the choice of the subdivision $0=t_{0}<\cdots<t_{n}=1$ nor on the covering $\{U_{j}\}.$ 

One can verify that two homotopic loops in $M$ give homotopic $W_{\Sigma}$-loops. This enables us to define the homomorphism $f: \pi_{1}(M,p)\rightarrow \pi_{1}(W_{\Sigma},p).$

Finally we have to prove that the homomorphism $f$ defined above is surjective.
It suffices  to show that we can lift a $W_{\Sigma}$-loop  $(w_{i},c_{i})_{0\leq i\leq n}$ in $\Sigma$ to a loop $\alpha$ in $M.$ By definition  $w_{i}=\varphi_{[\beta_{i}]},$ where $\beta_{i}$ is a curve contained in a regular leaf and  $\varphi_{[\beta_{i}]}$ is a parallel displacement map along $\beta.$
Define $\delta_{i}(t):=\varphi_{[\beta_{i}^{t}]}(c_{i}(t_{i})),$ where $\beta_{i}^{t}(s):=\beta_i(s\,  t)$ for $0\leq s\leq 1.$ It follows from Corollary \ref{cor-map-paralelo} that $\delta_{i}$ is a curve in the leaf $L_{c_{i}(t_{i})}.$ The concatenations of $\delta_{i}$ and $c_{i}$ give us the desired loop $\alpha.$
\end{proof}

\begin{cor}
\label{pi-M/F-trivial}
Suppose that the leaves of $\F$ are compact and the holonomy of each leaf is trivial. Let $\Sigma$ be a section of $\F$ and $p\in\Sigma.$  Then $\pi_{1}(M/\F)=\pi_{1}(W_{\Sigma},p).$ Furthermore, if $M$ is simply connected, so is $M/\F$.
\end{cor}

\begin{proof}
Proposition \ref{orbifoldcoxeter} implies that $M/\F=\Sigma/W_{\Sigma}$ is an orbifold. Now the result follows from the above proposition together with the fact that the fundamental group of the orbifold $\Sigma/W_{\Sigma}$ coincides with the fundamental groups of the pseudogroup $W_{\Sigma}$ (see Appendix).
\end{proof}

\begin{prop}
\label{prop-omega-cover-regular-space-of-leaf}
Assume that each regular leaf of $\F$  is compact and has trivial holonomy. Let $\Sigma$ be a section of $\F$ and $\Omega$ a connected component of $M_{r}\cap \Sigma.$ Then $\Pi:\Omega\rightarrow M_{r}/\F$ is a covering map.
\end{prop}
\begin{proof} Surjectivity of $\Pi:\Omega\rightarrow M_{r}/\F$ follows from the next lemma. Since the leaves of $\F$ are compact and have trivial normal holonomy, $\Pi:\Omega\rightarrow M_{r}/\F$ is a covering map. 
\end{proof}
\begin{lema}
\label{lemaconvex}
Assume that each regular leaf of $\F$  is compact and has trivial holonomy. Let $\Sigma$ be a section of $\F$ and $\Omega$ a connected component of $M_{r}\cap \Sigma.$
Let $p$ be a point of $\Omega,$ $L$ a regular leaf of $\F$ and $\gamma|_{[0,1]}$ a geodesic segment in $\Sigma$ such that $\gamma(0)=p,$ $\gamma(1)\in L$ and $\length(\gamma)=\dist(p,L).$ Then there are no singular points on $\gamma$.
\end{lema}

\begin{proof}
Suppose that there exists $0<s<1$ such that $\gamma(s)$ is a singular point. Since singular points are isolated on $\gamma$ (see Corollary \ref{pontos-singulares-na-geodesica}) there exists $\epsilon>0$ such that $\gamma(s+\epsilon)$ and $\gamma(s-\epsilon)$ are regular points. Let $P_{\gamma(s+\epsilon)}$ be the plaque  through $\gamma(s+\epsilon)$. It follows from the Slice Theorem that there exists $x\in S_{\gamma(s)}\cap P_{\gamma(s+\epsilon)}$ such that $\dist(\gamma(s-\epsilon),x)< \dist(\gamma(s-\epsilon),\gamma(s+\epsilon)).$ Indeed, we can define $x:=C\cap P_{\gamma(s+\epsilon)},$ where $C$ is the Weyl chamber that contains $\gamma(s-\epsilon)$ and hence the result follows from the convexity of $C.$ Let $\widehat{\gamma}_{1}$ be the minimal geodesic segment between $\gamma(s-\epsilon)$ and $x.$ 

Now let $\xi$ be a parallel normal field along the plaque $P_{\gamma(s+\epsilon)}$ such that  $\xi_{\gamma(s+\epsilon)}$ is multiple of  $\dot\gamma(s+\epsilon)$ with $\|\xi\|=\dist(\gamma(s+\epsilon),\gamma(1)).$ Define $\widehat{\gamma}_{2}(t):=\exp_{x}(t\xi).$ It follows from Corollary \ref{cor-map-paralelo}
 that $\widehat{\gamma}_{2}(1)\in L_{\gamma(1)}.$ Finally define $\gamma_{0}:=\gamma|_{[0,s-\epsilon]}.$  Now it is easy to see that the concatenation $ \gamma_{0}*\widehat{\gamma}_{1}*\widehat{\gamma}_{2}$ is shorter than $\gamma,$ contradiction.
\end{proof}

Let  $\Sigma$ be a section of $\F.$ Note that  the elements of $W_{\Sigma}$ send $\Sigma_{r}$ (the set of regular points in $\Sigma$) into $\Sigma_{r}.$ Therefore, we can define the pseudogroup $W_{\Sigma_{r}}$ as the elements of $W_{\Sigma}$ with domains and images in $\Sigma_{r}.$ 

\begin{prop}
\label{WloopsHomotopy}
Let $\Sigma$ be a section of $\F$ 
and let $p\in\Sigma.$ Consider two $W_{\Sigma}$-loops $\gamma_{0}$ and $\gamma_{1}$ based at $p$ that belong to the same homotopy class of $\pi_{1}(\WS, p).$ Suppose that $\gamma_{0}$ and $\gamma_{1}$ are contained in $\Sigma_{r}.$
Then $\gamma_{0}$ and $\gamma_{1}$ belong to the same homotopy class of $\pi_{1}(\WSr, p).$ 
\end{prop}

\begin{proof}
According to Corollary \ref{transnormalprojection}, there exists an open covering of $M$ by distinguished neighborhoods $\{U_{k}\}$ such that the plaques in each neighborhood $U_{k}$ are the preimages of a continuous map $T_{k}:U_{k}\rightarrow C_{k},$ where $C_{k}$ is a Weyl chamber in $U_{k}.$

Since $\gamma_{0}$ and $\gamma_{1}$ belong to the same homotopy class of $\pi_{1}(\WS,p),$ there exist partitions $0=s_{0}<\cdots<s_{n}=1,$ $0=t_{0}<\cdots<t_{n}=1,$ elements $w_{i,j}\in \WS$ and homotopies $h_{i,j}:[s_{i-1},s_{i}]\times [t_{j-1},t_{j}]\rightarrow \Sigma$ with the following properties:
\begin{enumerate}
\item The image of each homotopy $h_{i,j}$ is contained in a neighborhood $U_{k}.$
\item $(w_{i,j},h_{i,j}(s,\cdot))_{1\leq j\leq n}$ is a $\WS$-loop, for each $s$ with $s_{i-1}\leq s\leq s_{i}.$
\item The $\WS$-loops $(w_{i,j},h_{i,j}(s_{i},\cdot))_{1\leq j\leq n}$ and $(w_{i+1,j},h_{i+1,j}(s_{i},\cdot))_{1\leq j\leq n}$ are equivalent.
\item $\gamma_{0}$ is equivalent to the $\WS$-loop $(w_{1,j},h_{1,j}(0,\cdot))_{1\leq j\leq n}.$
\item $\gamma_{1}$ is equivalent to the $\WS$-loop $(w_{n,j},h_{n,j}(1,\cdot))_{1\leq j\leq n}$
\end{enumerate} 

Since the image of each homotopy $h_{i,j}$ is contained in a neighborhood $U_{k},$ we can suppose that the image of $h_{i,j}$ is contained in a Weyl chamber $C_{i,j}\subset\Sigma.$ Indeed, the image of $T_{k}\circ h_{i,j}$ is contained in a Weyl chamber $C_{k}.$ We can find a parallel displacement map $\varphi$ such that $C_{i,j}:=\varphi(C_{k})\subset\Sigma.$

We will construct homotopies $\widehat{h}_{i,j}:[s_{i-1},s_{i}]\times [t_{j-1},t_{j}]\to M_{r}\cap\Sigma$ whose edges are close or identical to the edges of $h_{i,j}$ and which have the following properties:
 
\begin{enumerate}
\item The image of each homotopy $\widehat{h}_{i,j}$ is contained in the interior $\interior(C_{i,j})$ of the Weyl Chamber $C_{i,j}$. 
\item $(w_{i,j},\widehat{h}_{i,j}(s,\cdot))_{1\leq j\leq n}$ is a $\WSr$-loop, for each $s_{i-1}\leq s\leq  s_{i}.$
\item The $\WSr$-loops $(w_{i,j},\widehat{h}_{i,j}(s_{i},\cdot))_{1\leq j\leq n}$ and $(w_{i+1,j},\widehat{h}_{i+1,j}(s_{i},\cdot))_{1\leq j\leq n}$ are equivalent.
\item $\gamma_{0}$ is equivalent to the $\WSr$-loop $(w_{1,j},\widehat{h}_{1,j}(0,\cdot))_{1\leq j\leq n}.$
\item $\gamma_{1}$ is equivalent to the $\WSr$-loop $(w_{n,j},\widehat{h}_{n,j}(1,\cdot))_{1\leq j\leq n}$.
\end{enumerate}

We start with the construction of the homotopy $\widehat{h}_{1,1}.$

The first step is to find curves $\widehat{h}_{1,1}(\cdot,0),$ $ \widehat{h}_{1,1}(s_{1},\cdot)$
and $\widehat{h}_{1,1}(\cdot,t_{1})$  close to the curves $h_{1,1}(\cdot,0),$ $h_{1,1}(s_{1},\cdot)$ and $h_{1,1}(\cdot,t_{1})$
such that:

\begin{itemize}
\item The new curves  are contained  $\interior (C_{1,1})$ (interior of $C_{1,1}$).
\item The concatenation of $h_{1,1}(0,\cdot)$ and the new curves with appropriate orientation is a loop.
\item $\widehat{h}_{1,1}(s_{1},\cdot)$ is equivalent to a curve in the  interior of $C_{2,1}.$
\item $w_{1,1}\widehat{h}_{1,1}(\cdot,t_{1})$ is contained in the interior of $C_{1,2}.$
\end{itemize}

The second step is to find a homotopy $\widehat{h}_{1,1}:[0,s_{1}]\times[0,t_{1}]\rightarrow \interior (C_{1,1})$ whose edges are the new curves constructed in the first step. This is possible for 
the Weyl chamber is simply connected.

Similarly we construct the homotopy $\widehat{h}_{1,2}:[0,s_1]\times [t_1,t_2]\to \interior (C_{1,2})$.

The first step is to find curves $\widehat{h}_{1,2}(\cdot,t_{1}),$ $ \widehat{h}_{1,2}(s_{1},\cdot)$ and  $\widehat{h}_{1,2}(\cdot,t_{2})$  close to the curves $h_{1,2}(\cdot,t_{1}),$ $ h_{1,2}(s_{1},\cdot)$ and  $h_{1,2}(\cdot,t_{2})$ 
such that:

\begin{itemize}
\item The new curves  are contained in $\interior (C_{1,2}).$ 
\item The concatenation of $\widehat{h}_{1,2}(0,\cdot)$ and  the new curves with appropriate orientation  is a loop.
\item $\widehat{h}_{1,2}(\cdot,t_{1})=w_{1,1} \widehat{h}_{1,1}(\cdot,t_{1}).$
\item $\widehat{h}_{1,2}(s_{1},\cdot)$ is equivalent to a curve in the interior of $C_{2,2}.$
\item $w_{1,3}\widehat{h}_{1,2}(\cdot,t_{2})$ lies in the interior of $C_{1,3}.$
\end{itemize}

As before the second step is to find a homotopy $\widehat{h}_{1,2}:[0,s_{1}]\times[t_{1},t_{2}]\rightarrow \interior (C_{1,2})$ whose edges are the new curves constructed in the first step.

Now the construction of the other homotopies $\widehat{h}_{i,j}$ is straightforward. 

\end{proof}

\begin{cor}
\label{prop-homotopy}
 Suppose that the leaves of $\F$ are compact and the holonomy of each  regular leaf is trivial.  Let $H:[0,1]\times[0,1]\rightarrow M/\F$ be a homotopy such that $H(0,\cdot)$ and $H(1,\cdot)$ are curves in $M_{r}/\F.$ Then there exists a homotopy $\widehat{H}:[0,1]\times[0,1]\rightarrow M_{r}/\F$ such that $\widehat{H}(0,\cdot)=H(0,\cdot)$ and $\widehat{H}(1,\cdot)=H(1,\cdot).$ Furthermore, if $H$ fixes endpoints, then $\hat{H}$ fixes endpoints.
\end{cor}

\textbf{Convention:}\emph{
For the rest of this section, let $\F$ be a s.r.f.s. on a simply connected complete riemannian manifold $M,$ such that each regular leaf is compact and has trivial holonomy. Moreover, let $\Sigma$ be a section of $\F$ and $\Omega$ a connected component of $M_{r}\cap \Sigma.$
}

\begin{cor}
\label{pi-Mreg/F-trivial} The fundamental group $\pi_{1}(M_{r}/\F)$ is trivial.
\end{cor}

\begin{proof}
This follows directly from Corollary \ref{pi-M/F-trivial}
and Corollary \ref{prop-homotopy}.
\end{proof}

\begin{cor}
\label{cor-Pi-omega-Mreg/F-bijective}
The map $\Pi|_{\Omega}:\Omega\rightarrow M_{r}/\F$ is a homeomorphism. Therefore a regular leaf meets $\Omega$ exactly once.
\end{cor}

\begin{proof}
This follows directly from Proposition \ref{prop-omega-cover-regular-space-of-leaf}
and Corollary \ref{pi-Mreg/F-trivial}.
\end{proof}

\begin{prop}
\label{propconvex}
Let $p,q$ be points in $\Omega.$ Then each minimal geodesic segment in $M$ between $p$ and $q$ is contained in $\Omega.$ This means that $\Omega$ is convex in $M$ and therefore in $\Sigma$.
\end{prop}

\begin{proof}

Let $\gamma:[0,1]\to M$ be a geodesic segment from $p$ to $L_q$ of length equal to  $\dist(p,L_{q}).$ Then $\gamma$ intersects $L_q$ orthogonally at $t=1$ and is therefore contained in the section $\Sigma$. Then it follows from Lemma \ref{lemaconvex} that $\gamma$ is  contained in $\Omega.$ On the other hand $L_q$ intersects $\Omega$ only in $q$ by Corollary \ref{cor-Pi-omega-Mreg/F-bijective}. Therefore $\gamma(1)=q$ and $d(p,q)=d(p,L_q)$. Hence $\gamma$ is a minimal geodesic segment between $p$ and $q.$ For another minimal geodesic segment $\beta$ between $p$ and $q,$ we have $\length(\beta)=\dist(p,q)=\dist(p,L_{q}).$ Therefore the above argumentation applies to show that $\beta$ is also contained in $\Omega.$ 
\end{proof}

\begin{prop}
\label{PropOmegaIsOrbifold}
The map $\Pi:\overline{\Omega}\rightarrow M/\F$ is a homeomorphism. Therefore each leaf meets $\overline{\Omega}$ exactly once.
\end{prop}

\begin{proof}
%Let $\xi$ be a parallel normal field along a regular leaf $L.$ Define $\eta_{\xi}:L\rightarrow %M$ as $\eta_{\xi}(x):=\exp_{x}(\xi).$

First we want to prove that $\Pi:\overline{\Omega}\rightarrow M/\F$ is surjective.

Let $q\in\Sigma$ be a singular point, $\sigma$ a local section centered at $q$ and let $C$ be a Weyl chamber of $\sigma$ that contains $q.$ Choose a regular point $x$ in $C$  and let $\xi_{x}\in T_{x}\Sigma$ denote the smallest vector with $\exp_{x}(\xi_{x})=q.$  The corresponding geodesic from $x$ to $q$ has no singular point up to $q$. Let $\xi$ be the parallel normal translation of $\xi_{x}$ along $L.$ Since $\Pi:\Omega\rightarrow M_r/\F$ is bijective by Proposition \ref{cor-Pi-omega-Mreg/F-bijective}, there exists $\tilde{x}\in\Omega\cap L_{x}.$ On the other hand $\eta_{\xi}(L_{x})=L_{q}$ by Corollary \ref{cor-map-paralelo}. Therefore  $\eta_{\xi}(\tilde{x})\in\partial \Omega\cap L_{q}.$ 

Now we want to prove that $\Pi:\overline{\Omega}\rightarrow M/\F$ is injective.

Suppose that there exist two points $y,z\in\partial\Omega$ such that $\Pi(y)=\Pi(z).$ Let $\epsilon$ be small enough such that $B_{\epsilon}(y)\cap B_{\epsilon}(z)=\emptyset.$ Let $x$ be a regular point such that $x\in B_{\epsilon}(y)\cap \Omega.$ Let  $\xi_{x}\in T_{x}\Sigma$ denote the vector such that $\exp_{x}(\xi_{x})=y.$ Since $\eta_{\xi}:L_{x}\rightarrow L_{y}$ is surjective there exists $\tilde{x}\in L_{x}$ such that $\eta_{\xi}(\tilde{x})=z.$ Due to the slice theorem we can suppose that $\tilde{x}\in\Omega\cap B_{\epsilon}(z).$ Since $\tilde{x}\in L_{x},$ we have $\Pi(x)=\Pi(\tilde{x})$ contradicting bijectivity of $\Pi:\Omega\rightarrow M_r/\F$ (see Corollary \ref{cor-Pi-omega-Mreg/F-bijective}).
 \end{proof}

\subsection{Proof of theorem \ref{teo-fundamental-domain}}
The propositions above and Theorem \ref{teo-excepcional-leaf}
 allow us  to prove Theorem \ref{teo-fundamental-domain}.
In fact item 1 of Theorem \ref{teo-fundamental-domain} follows from Proposition \ref{orbifoldcoxeter} and Corollary \ref{pi-M/F-trivial}. Item 2 of Theorem \ref{teo-fundamental-domain} follows from Corollary \ref{cor-Pi-omega-Mreg/F-bijective},
 Proposition \ref{propconvex} and Proposition \ref{PropOmegaIsOrbifold}.

%%%%%%%%%%%%%%%
 %%%%%%%%%%%%%%
\section{The transverse frame bundle of $\F.$}

In \cite{Molino} Molino associated an $O(k)$-principal bundle, the {\it orthogonal transverse frame bundle} to a regular riemannian foliation $(M,\F)$ of codimension $k$. A fiber of this bundle over a point $p$ in $M$ is defined as the set of orthonormal $k$-frames in $\nu_p L_p$, where $L_p$ is the leaf through $p$. In this section we generalize this notion for a s.r.f.s. of codimension $k$. Its restriction to the regular stratum $M_r$ will coincide with the orthogonal transverse frame bundle in the sense of Molino.

%This section can be omitted on a first reading by one interested in the proof of 

The reader, who is mainly interested in the  proof of Theorem \ref{teo-excepcional-leaf}, can omit this section and proceed with a summary of properties in the next section.

In order to motivate the definition of the tranverse frame bundle associated to a s.r.f.s., we present the following example.

\begin{ex}
\label{trivialex} 
The simplest example of a s.r.f.s. is the one of $M:=\mathbf{R}^{2}$ foliated by circles centered at the origin. We denote it by $\F$. The only singular leaf is the origin and the sections are the lines through the origin. Excising the singular leaf we obtain a regular riemannian foliation $\F_r$ of $M_{r}:=\mathbf{R}^{2}-\{(0,0)\}$. 
Let $E_{T}$ be the orthogonal transverse frame bundle (in the sense of Molino) associated to $\F_{r}.$ It is not difficult to see that $E_{T}=M^{1}_{r} \amalg M^{-1}_{r},$ where $M^{i}_{r}:=(\mathbf{R}^{2}-\{(0,0)\})\times \{i\}$  for $i=1,-1.$ We can identify $M^{1}_{r}$ (respectively $M^{-1}_{r}$) with the unit normal field outward (respectively inward) oriented.
 Set  $E= M^{1}\cup M^{-1},$ where $M^{i}:= M_{r}^{i}\cup (\{(0,0)\}\times \{i\}).$ We  will define $E$ as the  transverse frame bundle associated to $\F.$ It is obvious that the restriction of $E$ to $\pi^{-1}(M_r)$ is the orthogonal tranverse frame bundle $E_{T}.$ 

We would like to define a  structure $\tilde{q}$  associated to $q:=(0,0)$ that can be identified with the point $(0,0)\times 1$ of  $M^{1}.$ 
Since $M^{1}_{r}$ is identified with the outward oriented unit normal field, it is natural to look for a structure  that induces an outward orientation of $\F_{r}.$

If we start with a Weyl Chamber $[0,\infty)$ and the vector $(1,0)_{q}$ we can induce the desired orientation by parallel transport along the Weyl chamber and the circles.
On the other hand, if we start with the Weyl Chamber $(-\infty,0]$ and with a vector $(-1,0)_{q},$ we can induce the same orientation. This means that the pairs $((1,0)_{q},[0,\infty))$ and $((-1,0)_{q},(-\infty,0])$ represent the same object $\tilde{q}$, or more precisely, they belong to the same equivalence class.
This suggests the following definition. Let $(\zeta^{i}_{p},C_{i})$ for $i=1,2$ be a pair of a vector  tangential to some local sections $\sigma$ through $p$ with footpoint $p\in M$ and Weyl Chamber $C_{i}$  in $\sigma$ that contains $p.$ Then $(\zeta^{i}_{p},C_{i})$ are defined to be equivalent if there exists a rotation $\varphi$ (a parallel displacement map) such that $\varphi (C_{1})=\varphi (C_{2})$ and $\varphi_*\zeta^{1}=\varphi_*\zeta^{2}.$

We say that $\tilde{p}:=[(\zeta_{p},C)]$ belongs to $M^{1}$ if a representative $(\zeta_{p},C)$ induces the outward orientation of $\F_{r}.$
Note that if $p$ is not $(0,0)$ then there exists only one Weyl chamber $C$ that contains $p$ and hence this new definition coincides with the definition of a vector $\zeta_{p}$ with footpoint $p,$ when $p$ is regular.  

We would like to point out some particular and general aspects of this example.
In general, the transverse frame bundle $E$ associated to a s.r.f.s. $\F$ is not a union of the copies of $M.$, even if $\F$ is regular.
However, if the sections of a s.r.f.s. are flat,  there will be an integrable distribution on $E$ and hence $E$ will be foliated by leaves that cover $M$ (see Remark \ref{rem-srfsflatsections}).

The definition of $\tilde{p}$  will be valid  in the general case. The idea that $\tilde{p}$ induces frames, by parallel translation along the Weyl chambers and leaves,  will  play an important role. It will allow us to define  cross sections, which are useful, for example, to define the topology of $E.$  
\end{ex}
%%%%%%%%
\begin{dfn}[Tranverse Frame Bundle]
Let $\F$ be a s.r.f.s. on a complete riemannian manifold $M.$ 
Let $(\zeta_p,C)$ be a pair of an orthonormal $k$-frame $\zeta$ with footpoint $p$ tangential to a local section $\sigma$ and the germ of a Weyl chamber $C$ of $\sigma$ at $p$.
We identify $(\zeta^1_p,C_1)$ and $(\zeta^2_p,C_2)$ if there is a parallel displacement map $\varphi\in W_{S_p}$ (which fixes $p$) that maps  $ C_1$ to $C_2$ as germs in $p$ and $\zeta^1_p$ to $\zeta^2_p$ at first order. In other words, the equivalence class $[(\zeta_p,C)]$ consists of the $W_{S_p}$-orbit $(\varphi_*\zeta_p,\varphi(C)), \varphi\in W_{S_p}$, where $S_p$ is a slice through $p$. We call an equivalence class $[(\zeta_p,C)]$ {\it transverse frame}, and the set $E$ of transverse frames {\it transverse frame bundle}.
%\begin{rem}
%The next aim is to construct an $O(k)$-principal bundle $E=O_k(M,\F)\to M$ over $M$. First we construct $E|\bar C_w$ over the closures of the %chambers and glue them together on the common border of two chambers, obtaining $E|\sigma$. Similarly we obtain $E|S$ and then %$E|U$. Let 
%$$
%E|S:={\bigg(} \bigcup_{g\in\Gamma_p}E|(g\bar C_0){\bigg)}/\sim,
%$$ 
%where $\xi_1\sim \xi_2$ for $\xi_i\in O(g_i\bar C_0), g_i\in \Gamma_p$ with footprints $p_i$ if $p_1=p_2$, a section $\sigma$ such that %%$g_i(\bar C_0)$ are chambers of $\sigma$ and if there exists a $w\in W_{p_1}=W_{p_2}$ such that $w\xi_1=\xi_2$, where $W_{p_i}$ is the %subgroup of $W_\sigma$ fixing $\partial C_{w_1}\cap\partial C_{w_2}$.\par
%\end{rem}
Let $\pi:E\to M$ be the footpoint map. The fiber $F_q=\pi^{-1}(q)$ is equal to the set of transverse frames $[(\zeta_q,C)]$. There is a natural right action of $O(k)$ on $E$ by $[(\zeta_q,C)]\cdot g:=[(\zeta_q\cdot g,C)]$. This action is well-defined and simply transitive on the fiber. Note that in each equivalence class there is only one representative with a given  Weyl chamber.
\end{dfn}

Given a transverse frame $\tilde{q}=[(\zeta_{q},C)]$ over a point $q$ we want to find a neighborhood $U$ of $q$ in $M$ and a map  
$\s:U\rightarrow \pi^{-1}(U)$ such that $\s(q)=\tilde{q}$ and $\pi\circ\s(x)=x$ for $x\in U$. 
This will become a cross-section if we establish that $\pi:E\to M$ is a bundle. In fact we will use this map to define local trivializations for $\pi$.

First we define the cross-section $\s:S\to E|S$ over a slice $S$ through $q$. Let $x\in S$ be arbitrary. Then there is  $w\in W_S$ such that $x\in wC$. We define 
$$
\s(x)=[(({\overset{1}{\underset{0}\parallel}}\gamma_x)w\zeta,wC)],
$$
where $\gamma_x$ is the unique minimal geodesic segment from $q$ to $x$ and $({\overset{1}{\underset{0}\parallel}}\gamma_x)\zeta$ is the parallel transport of the frame $\zeta$ along $\gamma_{x}.$ This is independent of the choice of $w$.
 
Next we extend $\s|S$ to a cross-section on $U$. Take an arbitrary $y\in U$. Then there is $x\in S$ such that $y\in P_x$, where $P_x$ is the plaque of $U$ through $x$. Let $\varphi$ be a parallel displacement map along a curve in a regular plaque of $U$ such that $\varphi(x)=y$. Let $(\zeta_x,C)$ such that $\s(x)=[(\zeta_x,C)]$. We define $\s(y):=[(\varphi_*\zeta_x,\varphi(C))]$. 
Note that this expression is independent of the choice of $x$ and $\varphi.$ This follows from the definition of $\s$ on the slice $S_{q}$ and from the fact that the holonomy of each regular plaque in $U$ is trivial. 

The following map will become a trivialization of $E|U$. 
\begin{eqnarray}
\label{trivializationphi}
\phi:U\times O(k)&\to& E|U\nonumber\\
 (x,g)&\mapsto & \s(x)\cdot g.
\end{eqnarray}
Let $E|U$ take the induced topology via $\phi$. We have to show that this topology on $E$ is coherently defined, i.e., the transition from one trivialization to another trivialization is a homeomorphism. For this purpose it suffices to consider two cross-sections $\s_i:U_i\to E$ with $U_i\cap U_j\neq\emptyset$ and the corresponding trivializations $\phi_i:U_i\times O(k)\to E|U_i$. We define $h:U_1\cap U_2\to O(k)$ by $\s_1(x)=\s_2(x)\cdot h(x)$. 

\begin{prop}
\label{trivializationchange}
$\phi_2^{-1}\circ\phi_1(x,g)=(x,h(x)\cdot g)$.
\end{prop}
\begin{proof} For fixed $x\in U_1\cap U_2$ we write $\s_1(x)=[(\zeta^1_x,C)]$ and $\s_2(x)=[(\tilde\zeta_x,\tilde C)]$. There is $\zeta^2_x$ such $\s_2(x)=[(\zeta^2_x,C)]$. By assumption $\s_1(x)=\s_2(x)\cdot h(x)$, so $[(\zeta^1_x,C)]=[(\zeta^2_x \cdot h(x),C)]$. This implies
\begin{eqnarray}\label{formel1}
\zeta^1_x=\zeta^2_x\cdot h(x).
\end{eqnarray}
The definition of the trivialization $\phi$ (see Equation \ref{trivializationphi}) and that of the $O(k)$-action on a transverse frame together with Equation \ref{formel1} imply the next equation:
\begin{eqnarray*}
\phi_1(x,g)&=&\s_1(x)\cdot g\\
           &=&[(\zeta^1_x\cdot g,C)]\\
           &=&[(\zeta^2_{x} \cdot h(x) g,C)]\\
           &=&\s_2(x)\cdot h(x)g\\
	     &=&\phi_2(x,h(x)g)	
\end{eqnarray*}
This proves the claim $\phi_2^{-1}\circ\phi_1(x,g)=(x,h(x)g)$.
% Therefore 
%$[(\zeta^1_xg,C)]=[(\zeta^2h(x)g,C)]$. It follows %$\phi_1(x,g)=\s_1(x)g=\s_2(x)h(x)g=\phi_2(x,h(x)g)$. This proves the claim $\phi_2^{-%1}\phi_1(x,g)=(x,h(x)g)$.
\end{proof}

\begin{prop}
$h$ is constant along the plaques in $U_1\cap U_2$.
\end{prop} 
\begin{proof} 
For fixed $x$ let $\sigma$ be a local section centered at $x$ and $C$ a Weyl chamber of $\sigma$ containing $x$. We write 
$$
\s_1(x)=[(\zeta_x^1,C)]\quad \mbox{and}\quad \s_2(x)=[(\zeta_x^2,C)].
$$
%We have seen 
%\begin{eqnarray}\label{formel1}
%\zeta^1_x=\zeta^2_x\cdot h(x).
%\end{eqnarray}
Let $\varphi$ be a parallel displacement map along a curve in a regular plaque of $U_1\cap U_2$. We want to prove that $h(x)=h(\varphi(x))$. By definition of $\s_i$
$$
\s_i(\varphi(x))=[\varphi_*\zeta^i_x,\varphi(C)],
$$
by definition of $h$ we have
\begin{eqnarray}\label{formel2}
\varphi_*\zeta^1_x=\varphi_*\zeta^2_x\cdot h(\varphi(x)).
\end{eqnarray}
On the other hand by (\ref{formel1}) we have
\begin{eqnarray}\label{formel3}
\varphi_*\zeta^1_x=\varphi_*\zeta^2_x\cdot h(x).
\end{eqnarray}
Equations (\ref{formel2}) and (\ref{formel3}) imply $h(x)=h(\varphi(x))$.
\end{proof}

\begin{prop}
The map $h:U_1\cap U_2\to O(k)$ is continuous at all points and differentiable at all regular points. If the sections are flat, $h$ is locally constant.
\end{prop}
\begin{proof} We prove continuity. First we recall that in order to define the cross-section $\s_{i}$ we start with a tranverse frame $\tilde{q_{i}}=[(\zeta_{q_{i}}^{i},C_{i})]$ over some point $q_{i}\in M.$ 
Let $\sigma$ be a local section centered at $q_1$ with Weyl chamber $C_{1}$. Let $z$ be the unique point of intersection of $P_{q_{2}}$ with $C_{1}$. Let $\zeta^2_{z}$ be the  frame at $z$ such that $\s_2(z)=[(\zeta^2_{z},C_{1})]$. Let $e^1$ be the geodesic frame in $\sigma$ centered at $q_1$ with $e^1(q_1)=\zeta^1_{q_{1}},$ i.e., the frame in $\sigma$ defined by parallel translations of $ \zeta^1_{q_{1}}$ along radial geodesics starting in $q_{1}.$ Let $e^2$ be the geodesic frame centered at $z$ with $e^2(z)=\zeta^2_{z}$. We define $\widehat{h}$ on $\sigma$ by
$$
e^1(x)=e^2(x)\cdot\widehat{h}(x).
$$
Clearly $\widehat{h}$ is differentiable. If the sections are flat, then $\widehat{h}$ is constant.  Note that 
$[(e^i(x),C_{1})]=\s_i(x)$ for all $x\in C_{1}$. Thus $h(x)=\widehat{h}(x)$ for all $x\in C_{1}$. Since $h$ is constant along the leaves, $h$ is continuous.
\end{proof}

\begin{dfn}[Singular $C^0$-Foliation on $E$]
Next we define a singular partition $\widetilde{\F}$ on $E$ as follows: Let $\phi:U\times O(k)\to E|U$ be a trivialization and $P_x$ for $x\in U$ the plaque of $\F$ in $U$. We define $\widetilde{\F}|U$ by the partition $\widetilde{P}_{\phi(x,g)}:=\phi(P_x,g)$. Since the transition map $h$ is constant along the plaques, $\widetilde{\F}$ is well-defined on $E$. We define a leaf $\widetilde{L}$ through a point $x$ as the set of endpoints of continuous paths contained in plaques that start in $x$. The restriction of $\widetilde{\F}$ to the bundle $E_r=E|M_r$ over the regular stratum $M_r$ is the standard foliation described in \cite{Molino}. 
\end{dfn}

\begin{prop}[Parallel Displacement Maps of $\widetilde{\F}$]
\label{displacementoftildeF}
Let $\tilde{\beta}$ be a curve in a leaf of $\widetilde{\F}$ and $\beta:=\pi\circ\tilde{\beta}$ the projection of $\tilde{\beta}.$ Let $\varphi_{[\beta]}:V_{0}\rightarrow V_{1}$ be a parallel displacement map of $\F.$ Then there exists a map $\widetilde{\varphi}_{[\tilde{\beta}]}:\pi^{-1}(V_{0})\rightarrow \pi^{-1}(V_{1})$ such that
\begin{enumerate}
\item $\varphi_{[\beta]}\circ \pi=\pi\circ\widetilde{\varphi}_{[\tilde{\beta}]}.$
\item $\widetilde{\varphi}_{[\tilde{\beta}]}(x)\in \widetilde{L}_{x}.$
\end{enumerate}
\end{prop}
$\widetilde{\varphi}_{[\tilde{\beta}]}$ is called a parallel displacement map.

\begin{proof}
We can choose a partition $0=t_{0}<\cdots<t_{n}=1$ such that $\beta_{i}:=\beta|_{[t_{i-1},t_{i}]}$ is contained in a neighborhood $U_{i}$ for which there exists a trivialization $\phi: U_{i}\times O(k)\rightarrow \pi^{-1}(U_{i}).$ Set $(\pi(x),g):=\phi^{-1}(x)$ for $x\in\pi^{-1}(U_{i}).$ Then define $\widetilde{\varphi}_{[\tilde{\beta_{i}}]}(x):=\phi(\varphi_{[\beta_{i}]}(\pi(x)),g)$
It follows from Proposition \ref{trivializationchange} that $\widetilde{\varphi}_{[\tilde{\beta}_{i}]}$ is well-defined, i.e., it does not depend on the choice of the trivialization $\phi.$ Finally define $\widetilde{\varphi}_{[\tilde{\beta}]}:=\widetilde{\varphi}_{[\tilde{\beta}_{n}]}(x)\circ \cdots\circ\widetilde{\varphi}_{[\tilde{\beta}_{1}]}(x)$
\end{proof}

\begin{rem}[Weyl Pseudogroup $W$]
\label{defWeypseudogrouptildeF}
Let $\Sigma$ be a section and $\sigma$ a local section. Set $\widetilde{\Sigma}:=\pi^{-1}(\Sigma)$ and $\widetilde{\sigma}:=\pi^{-1}(\sigma).$ Then we can define the pseudogroups $\WSt$ and $W_{\widetilde{\sigma}}$ as in Definition \ref{definitionWeylPseudogroup}.
\end{rem}

\begin{prop}[Local Trivialization of $\widetilde{\F}$]
\label{trivializationtildeF}
Let $q$ be a point of $M.$ Then there exists a neighborhood $U$ of $q,$ a Weyl Chamber $C$ contained in $U$ and a continuous surjective map  $\widetilde{T}:\pi^{-1}(U)\rightarrow \pi^{-1}(C),$ such that  $\widetilde{T}^{-1}(p)$ is a plaque of $\widetilde{\F}$ contained in $\pi^{-1}(U)$ for each $p\in\pi^{-1}(C).$ 
\end{prop}
The map $\widetilde{T}$ is called a \emph{local trivialization of} $\widetilde{\F}.$
\begin{proof}
Let $\phi:U\times O(k)\rightarrow \pi^{-1}(U)$ be a trivialization of the bundle $E$ and $T:U\rightarrow C$ be a trivialization of $\F.$ Set $(\pi(x),g):=\phi^{-1}(x).$ Define $\widetilde{T}(x):=\phi(T(\pi(x)),g).$ Proposition \ref{trivializationchange} implies that $\widetilde{T}$  is well-defined.
\end{proof}

\begin{rem}[Compact Leaves]
\emph{Let $\F$ be a s.r.f.s. with compact leaves on a complete riemannian manifold $M.$ Then the leaves of $\widetilde{\F}$ are compact in the tranverse frame bundle $E.$} Indeed, it is known that the regular leaves of $\widetilde{\F}$ are compact if the regular leaves of $\F$ are compact (see Molino \cite{Molino}, Proposition 3.7, page 94). Therefore, it suffices to check that a singular leaf $\widetilde{L}_{\tilde{q}}$ is compact. In what follows we sketch the proof of this case.

Set $q:=\pi(\tilde{q})$ and $L_{q}:=\pi(\widetilde{L}_{\tilde{q}})$ where $\pi:E\rightarrow M$ is the projection of the fiber bundle. Let $\tub(L_{q})$ be a saturated tubular neighborhood of $L_{q}.$ 

Suppose that $\widetilde{L}_{\tilde{q}}$ is not compact. Then there exists a curve $\tilde{\beta}$ contained in $\widetilde{L}_{\tilde{q}}$ and a sequence $t_{i}\rightarrow \infty $ such that the points $\tilde{\beta}(t_{i})$ are different from each other and $\pi(\tilde{\beta}(t_{i}))=q.$ Sliding along the leaves of $\widetilde{\F}$ following $\tilde{\beta},$ we can find a curve $\tilde{\alpha}$ such that $\pi\circ\tilde{\alpha}$ is contained in a regular leaf $L_{p}\subset \tub(L_{q}).$ We can suppose that $L_{p}$ has trivial holonomy (see Theorem \ref{frss-eh-equifocal}). The fact that the leaf $L_{p}$ is compact implies that there exists a subsequence $(t_{n_{i}})$ such that $x=\pi(\tilde{\alpha}(t_{n_{i}}))=\pi(\tilde{\alpha}(t_{n_{j}}))$ for all $i,j.$ Hence the curve $\pi\circ\tilde{\alpha}$ gives us an infinite number of loops $\alpha_{n_{i}}$ in $L_{p}$ based at $x.$ Finally, since the points $\tilde{\beta}(t_{n_{i}})$ are different from each other,  the holonomy of the loops $\alpha_{n_{i}}$ are different, i.e., $L_{p}$ has nontrivial holonomy,  which is a contradiction.   

\end{rem}

\begin{rem}[S.R.F.S. with Flat Sections]
\label{rem-srfsflatsections}
If the sections of $\F$ are flat, $E$ is a smooth bundle and $\widetilde{\F}$ is a (smooth) singular foliation. We can define a distribution $H$ on $E$ by $H_{\tilde{q}}:=T_{\tilde{q}}\s(U).$ where $\s:U\rightarrow E$ is the cross-section defined above with respect to $\tilde q$.
It is not difficult to check that this distribution is integrable. This implies that $E$ is foliated by submanifolds $\{\widetilde{M}_{\tilde{x}}\}$ and for each $\tilde x\in E$ the map $\pi:\widetilde{M}_{\tilde{x}}\rightarrow M$ is a covering map. For each manifold $\widetilde{M}_{\tilde{x}}$ the lift of $\F$ along $\pi$ coincides with $\widetilde{\F}|\widetilde{M}_{\tilde{x}}.$ This is exactly what happens in Example \ref{trivialex}.
The covering map $\pi:\widetilde{M}_{\tilde{x}}\rightarrow M$ is a diffeomorphism if $M$ is simply connected. This implies that the regular leaves of $\F$ have trivial holonomy. Hence we have proved Theorem \ref{teo-excepcional-leaf} when $\F$ has flat sections. The general case will be proved in the next section.
\end{rem}

\section{Proof of Theorem \ref{teo-excepcional-leaf}}
In this section, we use the principal bundle $E,$ defined in the last section to prove Theorem \ref{teo-excepcional-leaf}. The reader, who skipped the last section, can read this section if he accepts the following facts.
\begin{enumerate}
\item There exists a continuous principal $O(k)$-bundle $E$ over $M$ that is associated to the s.r.f.s. $\F.$ The restriction of $E$ over $M_{r}$ (denoted by $E_{r}$) coincides with the usual orthogonal transverse frame bundle of the riemannian foliation $\F_{r}$(the restriction of $\F$ to $M_{r}$).
\item There exists a singular $C^0$-foliation $\widetilde{\F}$ on $E.$ The restriction of $\widetilde{\F}$ to $E_{r}$ coincides with the usual parallelizable foliation $\widetilde{\F}_{r}$ on $E_{r},$ which is a foliation with trivial holonomy whose leaves cover the leaves of $\F_{r}.$
\item There exist parallel displacement maps associated to $\widetilde{\F}$; hence we can define a Weyl pseudogroup $\WSt$ (see  Proposition \ref{displacementoftildeF}
 and Remark \ref{defWeypseudogrouptildeF}).
\item There exist local trivializations of $\widetilde{\F}$ (see  Proposition \ref{trivializationtildeF}).
\end{enumerate}

Let $L$ be a regular leaf, $p\in L$ and $\alpha$ a curve in $L$ such that $\alpha(0)=p=\alpha(1).$ 
Let $\zeta(t)$ be the parallel transport of an orthonormal frame $\zeta$ in $p$ along $\alpha$. Note that $\zeta(t)$  is contained in a regular leaf of the singular foliation $\widetilde{\F}$ in $E$. 

We want to show that  $\zeta(0)=\zeta(1).$ 

Since $M$ is simply connected we have a homotopy $G:[0,1]\times [0,1]\to M$ with
\begin{enumerate}
\item $G(0,t)=\alpha(t)$ for all $t\in [0,1]$.
\item $G(s,0)=G(s,1)=p$ for all $s$.
\item $G(1,t)=p$ for all $t$. 
\end{enumerate}
We define $\tilde p:=\zeta(0)$. Let $\pi:E\to M$ be the canonical projection of the transversal frame bundle $E$ of $\F$. We can lift $G$ to a homotopy $\tilde G:[0,1]\times [0,1]\to E$ with
\begin{enumerate}
\item $\tilde G(0,t)=\zeta(t)$ for all $t$.
\item $\tilde G(s,0)=\tilde{p}$ for all $s$.
\item $\tilde G(1,t)=\tilde{p}$ for all $t$.
\item $\pi\circ \tilde G(s,1)=p$ for all $s$.
\end{enumerate}
Let $\Sigma$ be the section of $\F$ that contains $p$ and define $\widetilde\Sigma:=\pi^{-1}(\Sigma)$.
Let $\rho:E\to E/\widetilde\F$ be the natural projection. 

Since the regular leaves of $\widetilde{\F}$ are compact and have trivial holonomy (see Molino \cite{Molino}, Proposition 3.7, page 94), we have the following claim: 
\begin{claim}
\label{claim1}
 $\rho:\widetilde\Sigma_r\to E_r/\tilde\F$ is a covering map, where  $\widetilde{\Sigma}_r=\widetilde{\Sigma}\cap E_r$.
\end{claim}

We define $R(s,t)=(1-t,1-s)$ and $\gamma_0:=\tilde G\circ R(0,\cdot)$. $\tilde G\circ R$ is a homotopy between $\gamma_0$ and the constant curve $\gamma_{\tilde p}\equiv \tilde p$. We have $\gamma_0(0)=\tilde p$ and $\pi\circ \gamma_0\equiv p.$ These facts imply  that $\gamma_0$ and $\gamma_{\tilde p}$ are contained in $\widetilde{\Sigma}_{r}.$ 

\begin{claim}
\label{claim1,5}
$\gamma_{0}$ is a $\WSrt$-loop based at $\tilde{p}$ and in particular a $\WSt$-loop based at $\tilde{p}.$
\end{claim}
\begin{proof}
$\gamma_{0}$ is a $\WSrt$-loop since $\gamma_{0}(1)=\zeta(1)\in \widetilde{L}_{\tilde{p}}.$
\end{proof}

\begin{claim}
\label{claim2}
$\gamma_{0}$ and the trivial $\WSt$-loop $\gamma_{\tilde p}$ belong to the same homotopy class of  $\pi_{1}(\WSt, \tilde{p}).$
\end{claim}
\begin{proof} Using the parallel displacement maps and trivializations of $\widetilde{\F},$ we can project the homotopy $\tilde G\circ R$ to $\WSt$ -loop deformations on $\widetilde{\Sigma}$ as we have done in Proposition \ref{PropHomOfFundamentalGroupofPseudogroup}.
\end{proof}

\begin{claim}
\label{claim3}
Consider two $\WSt$-loops $\delta_{0}$ and $\delta_{1}$ based at $\tilde{p}$ that belong to the same homotopy class of $\pi_{1}(\WSt, \tilde{p}).$ Suppose that $\delta_{0}$ and $\delta_{1}$ are contained in $\widetilde{\Sigma}_{r}.$
Then $\delta_{0}$ and $\delta_{1}$ belong to the same homotopy class of $\pi_{1}(\WSrt, \tilde{p}).$ 
\end{claim}
\begin{proof}
This is the same proof as for Proposition \ref{WloopsHomotopy}.
\end{proof}

Claims \ref{claim1,5}, \ref{claim2}, \ref{claim3} and the fact that $\pi_{1}(\WSrt,\tilde{p})=\pi_{1}(E_{r}/\widetilde{\F},\rho(\tilde{p}))$ 
imply that $\rho\circ\gamma_0$ and the constant curve $\rho\circ\gamma_{\tilde p}$ are homotopic in $E_r/\tilde\F$ fixing endpoints. The lift of this homotopy along the covering $\rho:\tilde\Sigma_r\to E_r/\F$ (see claim \ref{claim1}) to the curve $\gamma_0$ in $\tilde\Sigma_r$ is a homotopy to a constant curve fixing endpoints. Thus $\zeta(0)=\gamma_0(0)=\gamma_0(1)=\zeta(1)$.

\section{Appendix}

In this section we recall the definitions of pseudogroups, orbifolds, $W$-loops and the fundamental group of a pseudogroup. More details can be found in the Appendix D in \cite{Molino} written by E. Salem, on which our Appendix is based.

\begin{dfn}[Pseudogroup]
Let $\Sigma$ be a $C^{0}$ manifold. A \emph{pseudogroup} $W$ of transformations of $\Sigma$ is a collection of homeomorphisms $w:U\rightarrow W,$ where $U$ and $V$ are open subsets of $\Sigma$ such that:
\begin{enumerate}
\item If $w\in W$ then $w^{-1}\in W.$
\item If $w:U\rightarrow V$ and $\tilde{w}:\tilde{U}\rightarrow\tilde{V}$ belong to $W,$ then 
$\tilde{w}\circ w:w^{-1}(\tilde{U})\rightarrow\tilde{V}$ also belongs to $W.$
\item If $w:U\rightarrow V$ belongs to $W,$ then its restriction to each open subset $\tilde{U}\subset U$ also belongs to $W.$
\item If $w:U\rightarrow V$ is a homeomorphism between open subsets of $\Sigma$ which coincides in a neighborhood of each point of $U$ with an element of $W,$ then $w\in W.$
\end{enumerate} 
\end{dfn}

\begin{dfn}
Let $A$ be a family of local homeomorphisms of $\Sigma$ containing the identity map of $\Sigma.$ 
The pseudogroup obtained by taking the inverses of the elements of $A,$ the restriction to open sets of elements of $A,$ as well as their compositions and their unions, is called the {\it pseudogroup generated by $A.$}
\end{dfn}

An important example of a pseudogroup is the holonomy pseudogroup of a foliation, whose definition we now recall. Let $\F$ be a foliation of codimension $k$ on a manifold $M.$ Then $\F$ can be  described by an open cover $\{U_{i}\}$ of $M$ with submersions $f_{i}:U_{i}\rightarrow\sigma_{i}$ (where $\sigma_{i}$ is a submanifold of dimension $k$) such that there are diffeomorphisms $w_{i,j}:f_{i}(U_{i}\cap U_{j})\rightarrow f_{j}(U_{j}\cap U_{i})$ with $f_{j}=w_{i,j}\circ f_{i}.$ The elements $w_{i,j}$ acting on $\Sigma=\amalg \sigma_{i}$ generate a pseudogroup of transformations of $\Sigma$  called the \emph{holonomy pseudogroup of $\F$}.

In our work, we have a pseudogroup associated to a singular riemannian foliation with sections  and $\Sigma$ will be a fixed section (in particular it will be a connected submanifold).

\begin{dfn}[Orbifold]
One can define a $k$-dimensional orbifold as an equivalence class of pseudogroups $W$ of transformations on a manifold $\Sigma$ (dimension of $\Sigma$ is equal to $k$) verifying the following conditions:
\begin{enumerate}
\item The space of orbits $\Sigma/ W$ is Hausdorff.
\item For each $x\in\Sigma,$ there exists an open neighborhood $U$ of $x$ in $\Sigma$ such that the restriction of $W$ to $U$ is generated by a finite group of diffeomorphisms of $U.$
\end{enumerate}
\end{dfn}

An important example of an orbifold is the space of leaves $M/\F$ where $M$ is a riemannian manifold and $\F$ is a riemannian foliation on $M$ with compact leaves (see \cite{Molino}, Proposition 3.7, page 94). 

In Theorem \ref{teo-fundamental-domain}  we will prove  that the space of leaves $M/\F$ is an orbifold, when $\F$ is a s.r.f.s. with compact leaves. In addition $M/\F$ is a Coxeter orbifold, i.e., for each $p  \in M/\F$ there exists a  neighborhood of $p$ in $M/\F$ that is homeomorphic to a Weyl Chamber of a Coxeter group. 

\begin{dfn}[W-loop]
A $W$-loop with base point $x_{0}\in\Sigma$ is defined by    
\begin{enumerate}
\item a sequence $0=t_{0}<\cdots<t_{n}=1,$    
\item continuous paths $c_{i}:[t_{i-1},t_{i}]\rightarrow \Sigma,$ $1\leq i\leq n,$  
\item elements $w_{i}\in W$ defined in a neighborhood of $c_{i}(t_{i})$ for $1\leq i \leq n$ such that $c_{1}(0)=w_{n}c_{n}(1)=x_{0}$ and $w_{i}c_{i}(t_{i})=c_{i+1}(t_{i}),$ where $1\leq i\leq n-1.$
\end{enumerate}
A \emph{subdivision} of such a $W$-loop is obtained by adding new points to the interval $[0,1],$ by taking the restriction of the $c_{i}$ to these new intervals and $w=id$ at the new points.
\end{dfn}

We now want to define homotopy classes of $W$-loops.

\begin{dfn}
Two $W$-loops are \emph{equivalent} if there exists a subdivision common to the loops represented by $(w_{i},c_{i})$ and $(\tilde{w}_{i}, \tilde{c}_{i})$ and elements $g_{i}\in W$ 
defined in a neighborhood of the path $c_{i}$ such that
\begin{enumerate}
\item  $g_{i}\circ c_{i}=\tilde{c}_{i},$ $1\leq i\leq n$,
\item $\tilde{w}_{i}\circ g_{i}$ and $g_{i+1}\circ w_{i}$ have the same germ at $c_{i}(t_{i}),$ $1\leq i\leq n-1,$
\item $\tilde{w}_{n}\circ g_{n}$ has the same germ at $c_{n}(1)$ as $w_{n}.$
\end{enumerate}
\end{dfn}

\begin{dfn}
A \emph{deformation} of a $W$-loop represented by $(w_{i},c_{i})$ is given by continuous deformations $c_{i}(s,\cdot)$ of the paths $c_{i}=c_{i}^{0}:[t_{i-1},t_{i}]\rightarrow \Sigma,$ such that $(w_{i},c_{i}(s,\cdot))$ represents a $W-$loop.
\end{dfn}

\begin{dfn}[Fundamental Group of a Pseudogroup]
Two $W$-loops are in the same \emph{homotopy class} if one can be obtained from the other by a series of subdivisions, equivalences and deformations.
The homotopy classes of $W$-loops based at $x_{0}\in \Sigma$ form a group $\pi_{1}(W,x_{0})$ called \emph{fundamental group of the pseudogroup} $W$ at the point $x_{0}.$
\end{dfn}

\begin{rem}
If the orbit space $\Sigma/W$ is connected, then there exists an isomorphism, defined up to conjugation, between $\pi_{1}(W,x)$ and $\pi_{1}(W,y)$ for $x,y$ in $\Sigma.$ In addition, if $\Sigma/W$ is a connected orbifold, then $\pi_{1}(W,x)=\pi(\Sigma/W, \rho(x)),$ where $\rho:\Sigma\rightarrow \Sigma/W$ is the natural projection.
\end{rem}

% pag 266, (pseudogrupo)
%pag271 def Orbifold
%pag 273 274 grupo fundamental de um pseudogrupo
%pag 274 obs ii 

\bibliographystyle{amsplain}

\end{document}